\input amstex
\input amsppt.sty
\magnification=\magstep1
\hsize=32truecc
\vsize=22.2truecm
\baselineskip=16truept
\NoBlackBoxes
\TagsOnRight \pageno=1 \nologo
\def\Z{\Bbb Z}

\def\Q{\Bbb Q}
\def\R{\Bbb R}

\def\l{\left}
\def\r{\right}
\def\bg{\bigg}
\def\({\bg(}
\def\[{\bg\lfloor}
\def\){\bg)}
\def\]{\bg\rfloor}
\def\t{\text}
\def\f{\frac}

\def\p{\ (\roman{mod}\ p)}

\def\bi{\binom}
\def\eq{\equiv}

\def\ls{\leqslant}

\def\mo{\roman{mod}}
\def\sign{\roman{sign}}
\def\ord{\roman{ord}}

\def\ve{\varepsilon}

\def\Proof{\noindent{\it Proof}}

\def\Remark{\medskip\noindent{\it  Remark}}

\def\Ack{\medskip\noindent {\bf Acknowledgments}}
\hbox {Finite Fields Appl. 56(2019), 285--307.}
\bigskip
\topmatter
\title On some determinants with Legendre symbol entries\endtitle
\author Zhi-Wei Sun\endauthor
\leftheadtext{Zhi-Wei Sun} \rightheadtext{On some determinants with Legendre symbol entries}
\affil Department of Mathematics, Nanjing University\\
 Nanjing 210093, People's Republic of China
  \\  zwsun\@nju.edu.cn
  \\ {\tt http://math.nju.edu.cn/$\sim$zwsun}
\endaffil
\abstract In this paper we mainly focus on some determinants with Legendre symbol entries.
Let $p$ be an odd prime and let $(\frac{\cdot}p)$ be the Legendre symbol.
We show that
$(\f{-S(d,p)}p)=1$ for any $d\in\Z$ with $(\f dp)=1$, and that
$$\left(\f{W_p}p\right)=\cases(-1)^{|\{0<k<\f p4:\ (\f kp)=-1\}|}&\text{if}\ p\eq1\pmod4,
\\(-1)^{\lfloor(p+1)/8\rfloor}&\text{if}\ p\eq3\pmod4,\endcases$$
where $$S(d,p)=\det\l[\l(\f{i^2+dj^2}p\r)\r]_{1\ls i,j\ls(p-1)/2}$$
and $$W_p=\det\l[\l(\f{i^2-((p-1)/2)!j}p\r)\r]_{0\ls i,j\ls(p-1)/2}.$$
We also pose some conjectures on determinants, one of which states that $(-1)^{\lfloor(p+1)/8\rfloor}W_p$ is a square when $p\equiv 3\pmod4$.
\endabstract
\thanks 2010 {\it Mathematics Subject Classification}.\,Primary 11C20;
Secondary 15A15, 11A07, 11R11. \indent {\it Keywords}. Legendre symbols; determinants, congruences modulo primes, quadratic fields.
\newline \indent The initital version of this paper was posted to arXiv with the ID {\tt arXiv:1308.2900} in 2013. This work was supported by the National Natural Science
Foundation of China (grant no. 11571162) and the NSFC-RFBR Cooperation and Exchange Program (grant no. 11811530072).
\endthanks
\endtopmatter
\document

\heading{1. Introduction}\endheading

For an $n\times n$ matrix $A=(a_{ij})_{0\ls i,j\ls n-1}$ over the field of complex numbers, we often write
the determinant $\det A$ (or $|A|$)
in the form $|a_{ij}|_{0\ls i,j\ls n-1}$. In this paper we study determinants with Legendre symbol entries.

Let $p$ be an odd prime and let $(\f{\cdot}{p})$ be the Legendre
symbol. The circulant determinant
$$\bg|\l(\f{j-i}p\r)\bg|_{0\ls i,j\ls p-1}=\vmatrix (\f 0p)&(\f1p)&(\f2p)&\hdots&(\f{p-1}p)\\(\f{p-1}p)&(\f0p)&(\f1p)&\hdots&(\f{p-2}p)
\\\vdots&\vdots&\vdots&\ddots&\vdots\\(\f1p)&(\f 2p)&(\f3p)&\hdots&(\f0p)\endvmatrix$$
takes the value
$$\prod_{r=0}^{p-1}\sum_{k=0}^{p-1}\l(\f kp\r)(e^{2\pi ir/p})^k=0$$
since $\sum_{k=0}^{p-1}(\f kp)=0$. (See [K99, (2.41)] for the evaluation of a general circulant determinant.)
For the matrix $A=(a_{ij})_{1\ls i,j\ls p-1}$ with
$a_{ij}=(\f{i-j}p)$, L. Carlitz [C59, Theorem 4] proved that its characteristic
polynomial is
$$|xI_{p-1}-A|=\l(x^2-\l(\f{-1}p\r)p\r)^{(p-3)/2}\l(x^2-\l(\f{-1}p\r)\r),$$
where $I_{p-1}$ is the $(p-1)\times (p-1)$ identity matrix. Putting
$x=0$ in Carlitz's formula,  we obtain that
$$|A|=\l(-\l(\f{-1}p\r)\r)^{(p-1)/2}p^{(p-3)/2}=p^{(p-3)/2}.$$

For $m\in\Z$ let $\{m\}_p$ denote the
least nonnegative residue of an integer $m$ modulo $p$.
For any integer $a\not\eq0\ (\mo\ p)$, $\{aj\}_p\ (j=1,\ldots,p-1)$
is a permutation of $1,\ldots,p-1$, and its
sign is the Legendre symbol $(\f ap)$ by Zolotarev's lemma (cf.
[DH] and [Z]). Therefore, for any integer $d\not\eq0\ (\mo\ p)$ we have
$$\aligned\bg|\l(\f{i+dj}p\r)\bg|\Sb 0\ls i,j\ls p-1\endSb=&
\bg|\l(\f{i-\{-dj\}_p}p\r)\bg|\Sb 0\ls i,j\ls p-1\endSb
\\=&\l(\f{-d}p\r)\bg|\l(\f{i-j}p\r)\bg|\Sb 0\ls i,j\ls p-1\endSb
=0\endaligned\tag1.1$$
and
$$\bg|\l(\f{i+dj}p\r)\bg|\Sb 1\ls i,j\ls p-1\endSb=\l(\f{-d}p\r)\bg|\l(\f{i-j}p\r)\bg|\Sb 1\ls i,j\ls p-1\endSb
=\l(\f{-d}p\r)p^{(p-3)/2}.\tag1.2$$

Let $p$ be an odd prime. In 2004, R. Chapman [Ch04] used quadratic
Gauss sums to determine the values of
$$\bg|\l(\f{i+j-1}p\r)\bg|\Sb 1\ls i,j\ls
(p-1)/2\endSb\quad\t{and}\quad  \bg|\l(\f{i+j-1}p\r)\bg|\Sb 1\ls i,j\ls(p+1)/2\endSb.$$
Since $(p+1)/2-i+(p+1)/2-j-1\eq -(i+j)\ (\mo\ p)$, we see that
$$\bg|\l(\f{i+j-1}p\r)\bg|\Sb 1\ls i,j\ls(p-1)/2\endSb=\l(\f{-1}p\r)\bg|\l(\f{i+j}p\r)\bg|\Sb 1\ls i,j\ls(p-1)/2\endSb$$
and
$$\bg|\l(\f{i+j-1}p\r)\bg|\Sb 1\ls i,j\ls(p+1)/2\endSb
=\bg|\l(\f{i+j}p\r)\bg|\Sb 0\ls i,j\ls(p-1)/2\endSb.\tag1.3$$
Let $\ve_p$ and $h(p)$
stand for the fundamental unit and the class number of the real quadratic field $\Q(\sqrt p)$ respectively, and write
$$\ve_p^{(2-(\f 2p))h(p)}=r_p+s_p\sqrt{p}\ \t{with}\ r_p,s_p\in\Q.$$
Chapman [Ch12] conjectured that
$$\bg|\l(\f{j-i}p\r)\bg|\Sb 0\ls i,j\ls(p-1)/2\endSb=\cases
-r_p&\t{if}\ p\eq1\ (\mo\ 4),\\1&\t{if}\ p\eq3\ (\mo\ 4).\endcases\tag1.4$$
As Chapman could not solve this problem for several years, he called the determinant {\it evil} (cf. [Ch12]).
Chapman's conjecture on his ``evil" determinant was finally confirmed by M. Vsemirnov [V12, V13]
via matrix decomposition. We guess that if $p\eq3\pmod4$ then
$|x+(\f{j-i}p)|_{1\ls i,j\ls(p-1)/2}=x.$

Let $p\eq1\ (\mo\ 4)$ be a prime. In an unpublished manuscript [Ch03] written in 2003, Chapman conjectured that
$$s_p=\l(\f 2p\r)\bg|\l(\f{j-i}p\r)\bg|\Sb 1\ls i,j\ls (p-1)/2\endSb.\tag1.5$$
Note that (1.4) and (1.5) together yield an interesting identity
$$\ve_p^{(2-(\f 2p))h(p)}=\l(\f 2p\r)\bg|\l(\f{j-i}p\r)\bg|_{1\ls i,j\ls(p-1)/2}\sqrt p-\bg|\l(\f{j-i}p\r)\bg|_{0\ls i,j\ls(p-1)/2}.$$
Taking the norm with respect to the field extension $\Q(\sqrt p)/\Q$, we are led to the identity
$$\bg|\l(\f{j-i}p\r)\bg|_{0\ls i,j\ls(p-1)/2}^2-p\ \bg|\l(\f{j-i}p\r)\bg|_{1\ls i,j\ls(p-1)/2}^2=-1$$
since $N(\ve_p)=-1$ and $2\nmid h(p)$ (cf. Theorem 3 of [Co62, p.\,185] and Theorem 6 of [Co62, p.\,187]).
This provides an {\it explicit} solution to the Diophantine equation $x^2-py^2=-1$ with $x,y\in\Z$.

Now we state our first theorem on evaluations of certain determinants over the finite field $\Bbb F_p=\Z/p\Z$ with $p$ prime.

\proclaim{Theorem 1.1} Let $p$ be an odd prime. For $d\in\Z$ define
$$R(d,p):=\bg|\l(\f{i+dj}p\r)\bg|\Sb 0\ls i,j\ls(p-1)/2\endSb.\tag1.6$$
If $p\eq1\ (\mo\ 4)$, then
$$R(d,p)\eq\l(\l(\f dp\r)d\r)^{(p-1)/4}\ \f{p-1}2!\pmod p.\tag1.7$$
When $p\eq3\ (\mo\ 4)$, we have
$$R(d,p)\eq\cases(\f2p)\ (\mo\ p)&\t{if}\ (\f dp)=1,
\\1\ (\mo\ p)&\t{if}\ (\f dp)=-1.\endcases\tag1.8$$
For the general odd prime $p$, we have
$$R(-d,p)\eq\l(\f 2p\r)R(d,p)\pmod p,\tag1.9$$
and
$$\bg|\l(\f{i+dj+c}p\r)\bg|\Sb 0\ls i,j\ls(p-1)/2\endSb\eq R(d,p)\pmod p\quad\t{for all}\ c\in\Z.\tag1.10$$
\endproclaim

\Remark\ 1.1. (a) By Wilson's theorem, for any odd prime $p$ we have
$$\l(\f{p-1}2!\r)^2\eq (-1)^{(p-1)/2}\prod_{k=1}^{(p-1)/2}k(p-k)\eq (-1)^{(p+1)/2}\pmod p.$$

(b) For a general integer $d$, we have not found any interesting pattern for the exact values of $R(d,p)$ with $p>2$ prime. For example, below are the values of $R(4,p)$ with $p$ prime and $2<p\ls41$.
$$\gather R(4,3)=-1,\ R(4,5)=-2,\ R(4,7)=1,\ R(4,11)=32,\ R(4,13)=8,
\\ R(4,17)=-72,\ R(4,19)=-2^7\cdot3\cdot5,\ R(4,23)=-2^5\cdot2719,
\\ R(4,29)=-2^7\cdot7\cdot1301,\ R(4,31)=-2^7\cdot83\cdot1433,
 \\R(4,37)=2^9\cdot3\cdot14479,\ R(4,41)=2^9\cdot38913779.
\endgather$$
\medskip

\proclaim{Corollary 1.1} Let $p\eq1\ (\mo\ 4)$ be a prime, and write $\ve_p^{h(p)}=a_p+b_p\sqrt p$ with $2a_p,2b_p\in\Z$.
Then
$$a_p\eq-\f{p-1}2!\pmod p. \tag1.11$$
\endproclaim
\Proof. By (1.7) we have
$$R(1,p)\eq\f{p-1}2!\pmod p.$$
On the other hand, in view of (1.3) and Chapman [Ch04, Corollary 3],
$$R(1,p)=\bg|\l(\f{i+j-1}p\r)\bg|\Sb 1\ls i,j\ls(p+1)/2\endSb
=-\l(\f2p\r)2^{(p-1)/2}a_p.$$
So (1.11) is valid. \qed
\medskip

It is well known that for any odd prime $p$ the $(p-1)/2$ squares
$$1^2,2^2,\ldots,\l(\f{p-1}2\r)^2$$
give all the $(p-1)/2$ quadratic residues modulo $p$. So we think that it's natural to consider
some determinants with Legendre symbol (or Jacobi symbol) entries related to binary quadratic forms.

For any integer $d$ and odd integer $n>1$, we define
$$S(d,n):=\bg|\l(\f{i^2+dj^2}n\r)\bg|\Sb 1\ls i,j\ls(n-1)/2\endSb\tag1.12$$
and
$$T(d,n):=\bg|\l(\f{i^2+dj^2}n\r)\bg|\Sb 0\ls i,j\ls(n-1)/2\endSb,\tag1.13$$
where $(\frac{\cdot}n)$ is the Jacobi symbol.

\proclaim{Theorem 1.2} {\rm (i)} Let $n>1$ be an odd composite number. Then, for any $d\in\Z$ we have
$S(d,n)=T(d,n)=0$.

{\rm (ii)} Let $p$ be an odd prime and let $c,d\in\Z$ with $p\nmid c$. Then
$$S(c^2d,p)=\l(\f cp\r)^{(p+1)/2}S(d,p)\ \ \t{and}\ \ T(c^2d,p)=\l(\f cp\r)^{(p+1)/2}T(d,p).\tag1.14$$
Also,
$$\l(\f dp\r)=-1\ \ \t{implies}\ \ S(d,p)=0.\tag1.15$$
When $p\eq 1\pmod4$, we have
$$S(-d,p)=\l(\f 2p\r)S(d,p)\quad\t{and}\quad T(-d,p)=\l(\f 2p\r)T(d,p).\tag1.16$$

{\rm (iii)} Let $p$ be an odd prime and let $d\in\Z$. Then
$$\l(\f{T(d,p)}p\r)=\cases (\f 2p)&\t{if}\ (\f dp)=1,
\\1&\t{if}\ (\f dp)=-1.\endcases\tag1.17$$
Also,
$$T(-d,p)\eq\l(\f 2p\r)T(d,p)\pmod p\tag1.18$$
and
$$\bg|\l(\f{i^2+dj^2+c}p\r)\bg|\Sb 0\ls i,j\ls(p-1)/2\endSb\eq T(d,p)\pmod p\quad\t{for all}\ c\in\Z.\tag1.19$$
Moreover, when $(\f dp)=1$, we have
$$S(d,p)=\f 2{p-1}T(d,p)\ \ \t{and}\ \ \l(\f{-S(d,p)}p\r)=1.\tag1.20$$
\endproclaim
\medskip
\Remark\ 1.2. By Theorem 1.2(iii), for any odd prime $p$ the number $-S(1,p)=2T(1,p)/(1-p)$
is a quadratic residue mod $p$. We note that $-S(1,p)$ is not always a square, and will pose a related conjecture
in Section 4.
\medskip

{\it Example} 1.1. Note that
$$S(1,11)=\vmatrix -1&1&-1&-1&1\\1&-1&-1&1&-1\\-1&-1&-1&1&1\\-1&1&1&-1&-1\\1&-1&1&-1&-1\endvmatrix=-4^2\ \ \t{with}\ \l(\f{-S(1,11)}{11}\r)=1,$$
and
$$S(2,13)=\vmatrix 1&1&-1&-1&1&-1\\-1&1&1&1&-1&-1\\-1&1&1&-1&-1&1\\-1&-1&-1&1&1&1\\1&-1&1&-1&1&-1\\1&-1&-1&1&-1&1\endvmatrix=0
\ \ \t{with}\ \l(\f 2{13}\r)=-1.$$
\medskip

For any positive odd integer $n$ and integers $c$ and $d$, we define
$$(c,d)_n:=\bigg|\l(\f{i^2+cij+dj^2}n\r)\bigg|_{1\ls i,j\ls n-1}\tag1.21$$
and
$$[c,d]_n:=\bigg|\l(\f{i^2+cij+dj^2}n\r)\bigg|_{0\ls i,j\ls n-1}.\tag1.22$$
It is easy to see that
$(-c,d)_n=(\f{-1}n)(c,d)_n$.

Now we present our third theorem.

\proclaim{Theorem 1.3} Let $c,d\in\Z$.

{\rm (i)} We have $(c,d)_n=0$ for any positive odd integer $n$ with $(\f dn)=-1$.

{\rm (ii)} If $p$ is an odd prime with $(\f dp)=1$, then
$$[c,d]_p=\cases\f{p-1}2(c,d)_p&\t{if}\ p\nmid c^2-4d,
\\\f{1-p}{p-2}(c,d)_p&\t{if}\ p\mid c^2-4d.\endcases\tag1.23$$
\endproclaim

Now we state two more theorems.

\proclaim{Theorem 1.4} {\rm (i)} For any odd prime $p$, we have
$$\bg|\f{(\f{i+j}p)}{i+j}\bg|_{1\ls i,j\ls (p-1)/2}\eq\cases(\f 2p)\ (\mo\ p)&\t{if}\ p\eq1\ (\mo\ 4),
\\((p-1)/2)!\ (\mo\ p)&\t{if}\ p\eq3\ (\mo\ 4).\endcases\tag1.24$$

{\rm (ii)} Let $p\eq3\ (\mo\ 4)$ be a prime. Then
$$\bg|\f1{i^2+j^2}\bg|_{1\ls i,j\ls (p-1)/2}\eq\l(\f 2p\r)\pmod p.\tag1.25$$
\endproclaim
\Remark\ 1.3. Similar to Theorem 1.4(ii), for any prime $p\eq5\ (\mo\ 6)$ we conjecture that
 $2|\f1{i^2-ij+j^2}|_{1\ls i,j\ls p-1}$
is a quadratic residue modulo $p$ and the $p$-adic order of $|\f1{i^2-ij+j^2}|_{1\ls i,j\ls (p-1)/2}$ is $(p+1)/6$. This has been verified for $p<500$.
\medskip

\proclaim{Theorem 1.5} For any odd prime $p$, we have
$$\left(\f{W_p}p\right)=\cases(-1)^{|\{0<k<\f p4:\ (\f kp)=-1\}|}&\text{if}\ p\eq1\pmod4,
\\(-1)^{\lfloor(p+1)/8\rfloor}&\text{if}\ p\eq3\pmod4,\endcases\tag1.26$$
where
$$W_p:=\left|\l(\f{i^2-\f{p-1}2!j}p\r)\r|_{0\ls i,j\ls (p-1)/2}.\tag1.27$$
\endproclaim
\Remark\ 1.4.  Below are the values of $W_p$ for the first 6 odd primes $p$:
$$W_3=1,\ W_5=-1,\ W_7=-4,\ W_{11}=-16,\ W_{13}=-288,\ W_{17}=-840.$$
See also Conjecture 4.1 inspired by Theorem 1.5.
\medskip

We are going to show Theorems 1.1-1.3 and Theorems 1.4-1.5 in Sections 2 and 3 respectively, and pose some new conjectures on determinants in Section 4.

\heading{2. Proof of Theorems 1.1-1.3}\endheading

\proclaim{Lemma 2.1 {\rm ([K05, Lemma 9])}} Let $P(z)=\sum_{k=0}^{n-1}a_kz^k$ be a polynomial with complex coefficients.
Then we have
$$|P(x_i+y_j)|_{1\ls i,j\ls n}=a_{n-1}^n\prod_{k=0}^{n-1}\bi{n-1}k\times\prod_{1\ls i<j\ls n}(x_i-x_j)(y_j-y_i).\tag2.1$$
\endproclaim
\Remark\ 2.1. This lemma remains valid over any commutative ring as pointed out by one of the referees.

\medskip
\noindent{\it Proof of Theorem 1.1}. Set $n=(p-1)/2$. For any $c\in\Z$, we have
$$\bg|\l(\f{i+dj+c}p\r)\bg|\Sb 0\ls i,j\ls(p-1)/2\endSb\eq|(i+dj+c)^n|\Sb 0\ls i,j\ls n\endSb\pmod p,$$
and
$$\align&|(i+dj+c)^n|\Sb 0\ls i,j\ls n\endSb=|(i+dj+c-d-1)^n|_{1\ls i,j\ls n+1}
\\=&\prod_{k=0}^n\bi nk\times\prod_{1\ls i<j\ls n+1}(i-j)(dj+c-d-1-(di+c-d-1))
\\=&\f{(n!)^{n+1}}{\prod_{k=0}^nk!(n-k)!}(-d)^{n(n+1)/2}\prod_{1\ls i<j\ls n+1}(j-i)^2
=(-d)^{n(n+1)/2}(n!)^{n+1}
\endalign$$
with the help of Lemma 2.1.
Therefore (1.10) holds, and also
$$R(d,p)\eq (-d)^{(p^2-1)/8}\l(\f{p-1}2!\r)^{(p+1)/2}\pmod p.\tag2.2$$
In the case $p\eq1\ (\mo\ 4)$, by (2.2) and Remark 1.1(a) we have (1.7) since
$$\align R(d,p)\eq& (-d)^{(p-1)/4\times (p+1)/2}\ \f{p-1}2!\l(\f{p-1}2!\r)^{2(p-1)/4}
\\\eq& \l(d^{(p+1)/2}\r)^{(p-1)/4}\ \f{p-1}2!\eq\l(\l(\f dp\r)d\r)^{(p-1)/4}\ \f{p-1}2!\pmod p.
\endalign$$
In the case $p\eq3\ (\mo\ 4)$, by (2.2) and Remark 1.1(a) we obtain
$$R(d,p)\eq (-d)^{(p-1)/2\times(p+1)/4}\l(\f{p-1}2!\r)^{2(p+1)/4}\eq\l(\f{-d}p\r)^{(p+1)/4}\pmod p$$
and hence (1.8) follows.

Now it remains to show (1.9) which holds trivially when $p\mid d$. If $p\eq 1\ (\mo\ 4)$, then by (1.7) we get
$$R(-d,p)\eq\(\l(\f{-d}p\r)(-d)\)^{(p-1)/4}\ \f{p-1}2!\eq\l(\f 2p\r)R(d,p)\pmod p.$$
If $p\eq3\ (\mo\ 4)$, then $(\f{-d}p)=-(\f dp)$ and hence we get (1.9) from (1.8).

The proof of Theorem 1.1 is now complete. \qed

Motivated by Zolotarev's lemma, in 2006 H. Pan obtained the following lemma.

\proclaim{Lemma 2.2 {\rm (H. Pan [P06])}} Let $n>1$ be an odd integer and let $c$ be any integer relatively prime to $n$. For each $j=1,\ldots,(n-1)/2$ let $\pi_c(j)$ be the unique $r\in\{1,\ldots,(n-1)/2\}$ with $cj$
congruent to $r$ or $-r$ modulo $n$. For the permutation $\pi_c$ on $\{1,\ldots,(n-1)/2\}$, its sign is given by
$$\sign(\pi_c)=\l(\f cn\r)^{(n+1)/2}.\tag2.3$$
\endproclaim

\proclaim{Lemma 2.3} Let $p\eq1\ (\mo\ 4)$ be a prime. Then
$$\l(\f{((p-1)/2)!}p\r)=\l(\f 2p\r).\tag2.4$$
\endproclaim
\Proof. Since
$$(-4)^{(p-1)/4}=(-1)^{(p-1)/4}\times2^{(p-1)/2}=\l(\f 2p\r)2^{(p-1)/2}\eq1\pmod p,$$
by Remark 1.1(a) there is an integer $x$ such that
$$x^4\eq-4\eq4\l(\f{p-1}2!\r)^2\pmod p,\ \ \t{i.e.,}\ x^2\eq\pm2\times\f{p-1}2!\pmod p.$$
Therefore (2.4) holds. \qed

\medskip

\noindent{\it Proof of Theorem} 1.2. (i) If $n=p^2$ for some odd prime $p$, then $1<p+1\ls (p+1)(p-1)/2=(n-1)/2$ and
$$\l(\f{1^2+dj^2}{p^2}\r)=\l(\f{(p+1)^2+dj^2}{p^2}\r)\quad\t{for all}\ j\in\Z.$$
If $n=pm$ with $p$ an odd prime and $m$ an odd integer greater than $p$, then for $i_1=(m-p)/2$ and $i_2=(m+p)/2$
we have
$$1\ls i_1<i_2<m=\f np<\f n2,\  \ \ i_2^2-i_1^2=pm=n,$$
 and
$$\l(\f{i_1^2+dj^2}{n}\r)=\l(\f{i_2^2+dj^2}{n}\r)\quad\t{for all}\ j\in\Z.$$
Therefore $S(d,n)=T(d,n)=0$.

(ii) Define the permutation $\pi_c$ as in Lemma 2.2 with $n=p$.
In view of (2.3), we have
$$S(c^2d,p)= \bg|\l(\f{i^2+d\pi_c(j)^2}p\r)\bg|_{1\ls i,j\ls {p-1}/2}=\l(\f cp\r)^{(p+1)/2}S(d,p)$$
as well as the second equality in (1.14).
Similarly, if $p\nmid d$ then
$$\l(\f dp\r)^{(p-1)/2}S(d,p)=\bigg|\l(\f{\pi_d(i)^2+dj^2}p\r)\bigg|_{1\ls i,j\ls(p-1)/2}=\l(\f dp\r)^{(p+1)/2}S(d,p)$$
and hence (1.15) holds.

Now we assume that $p\eq1\ (\mo\ 4)$. As $((p-1)/2)!^2\eq-1\pmod p$ by Remark 1.1(a),
applying (1.14) with $c=((p-1)/2)!$ and using (2.4) we immediately get (1.16).

(iii) For convenience we set $m=(p-1)/2$.
In light of Lemma 2.1 with $n=m+1$ and $P(z)=z^m$, for any $c\in\Z$ we have
$$\align&|(i^2+dj^2+c)^m|_{0\ls i,j\ls m}=|((i-1)^2+d(j-1)^2+c)^m|_{1\ls i,j\ls m+1}
\\=&\prod_{k=0}^m\bi mk\times\prod_{1\ls i<j\ls m+1}((i-1)^2-(j-1)^2)(d(j-1)^2+c-d(i-1)^2-c)
\\=&\f{(m!)^{m+1}}{\prod_{k=0}^mk!(m-k)!}(-d)^{m(m+1)/2}\prod_{0\ls i<j\ls m}(j-i)^2(j+i)^2
\\=&(-d)^{m(m+1)/2}(m!)^{m+1}\prod_{0\ls i<j\ls m}(i+j)^2.
\endalign$$
Therefore (1.19) holds, and also
$$\align T(d,p)\eq& (-d)^{(p^2-1)/8}\l(\f{p-1}2!\r)^{(p+1)/2}\prod_{0\ls i<j\ls (p-1)/2}(i+j)^2
\\\eq& R(d,p)\prod_{0\ls i<j\ls (p-1)/2}(i+j)^2\pmod p
\endalign$$
with the help of (2.2). Combining this with (1.9) we obtain (1.18). Note that
$$\l(\f{T(d,p)}p\r)=\l(\f{R(d,p)}p\r).\tag2.5$$
If $(\f dp)=1$, then by Theorem 1.1 we have
$$R(d,p)\eq\cases d^{(p-1)/4}((p-1)/2)!\ (\mo\ p)&\t{if}\ p\eq1\ (\mo\ 4),
\\(\f 2p)\ (\mo\ p)&\t{if}\ p\eq3\ (\mo\ 4),\endcases$$
and hence $(\f{R(d,p)}p)=(\f 2p)$ in view of Lemma 2.3.
In the case $(\f dp)=-1$, by Theorem 1.1 we have
$$R(d,p)\eq\cases (-d)^{(p-1)/4}((p-1)/2)!\ (\mo\ p)&\t{if}\ p\eq1\ (\mo\ 4),
\\1\ (\mo\ p)&\t{if}\ p\eq3\ (\mo\ 4),\endcases$$
and hence $(\f{R(d,p)}p)=1$ with the help of Lemma 2.3.
Therefore (1.17) also holds in light of (2.5).

Assume that $p\nmid d$. Then the sum of entries in each column of the determinant $S(d,p)$
actually equals $-(1+(\f dp))/2$.
In fact, for any $j=1,\ldots,m$, clearly
$$\bg|\sum_{i=1}^{m}\l(\f{i^2+dj^2}p\r)\bg|\ls m<\f p2$$
and also
$$\align\sum_{i=1}^{m}\l(\f{i^2+dj^2}p\r)\eq&\sum_{i=1}^{m}(i^2+dj^2)^{m}
\\=&\sum_{i=1}^m\(i^{2m}+d^mj^{2m}+\sum_{0<k<m}\bi mki^{2k}(dj^2)^{m-k}\)
\\\eq&m\l(1+\l(\f dp\r)\r)+\f12\sum_{0<k<m}\bi{m}k\sum_{i=1}^m (i^{2k}+(p-i)^{2k})(dj^2)^{m-k}
\\\eq&-\f{1+(\f dp)}2+\f12\sum_{0<k<m}\bi{m}k(dj^2)^{m-k}\sum_{i=1}^{p-1}i^{2k}
\\\eq&-\f{1+(\f dp)}2\pmod p.\endalign$$
(It is well known that $\sum_{i=1}^{p-1}i^h\eq0\pmod p$ for any $h\in\Z$ with $h\not\eq0\pmod{p-1}$,
see, e.g. [IR, p.\,235].)

Now suppose that $(\f dp)=1$. By the last paragraph, adding the last $m$ rows of $T(d,p)$ to the first row we see that
the initial term in the first row becomes $m$ while all the other terms in the first row turn out to be zero. It follows that
$$T(d,p)=m S(d,p)\eq-\f12S(d,p)\pmod p.$$
Combining this with $(1.17)$, we obtain that
$$\l(\f{-S(d,p)}p\r)=\l(\f{2T(d,p)}p\r)=1.$$
This proves (1.20). \qed

\medskip
\noindent{\it Proof of Theorem 1.3}. (i) Let $n$ be any positive odd integer relatively prime to $d$.
For $j=0,\ldots,n-1$ let $\sigma_d(j)$
be the least nonnegative residue of $dj$ modulo $n$. Then $\sigma_d$ is a permutation on $\{0,\ldots,n-1\}$ with $\sigma_d(0)=0$.
By Frobenius' extension of the Zolotarev lemma (cf. [BC15]), we have $\sign(\sigma_d)=(\f dn)$.

Now suppose that $(\f dn)=-1$. Then $n>1$ and $\sign(\sigma_d)=-1$. Thus
$$\align (c,d)_n=&\l(\f dn\r)^{n-1}(c,d)_n
=\bigg|\l(\f{di^2+ci(dj)+(dj)^2}n\r)\bigg|_{1\ls i,j\ls n-1}
\\=&\bigg|\l(\f{di^2+ci\sigma_d(j)+\sigma_d(j)^2}n\r)\bigg|_{1\ls i,j\ls n-1}
\\=&\sum_{\pi\in S_{n-1}}\sign(\pi)\prod_{i=1}^{n-1}\l(\f{di^2+ci\sigma_d(\pi(i))+\sigma_d(\pi(i))^2}n\r)
\\=&\sign(\sigma_d)\sum_{\tau\in S_{n-1}}\sign(\tau)\prod_{i=1}^{n-1}\l(\f{di^2+ci\tau(i)+\tau(i)^2}n\r)=-(c,d)_n
\endalign$$
and hence $(c,d)_n=0$.

(ii) Let $p$ be an odd prime. For any $a,b\in\Z$, it is known (cf. [BEW, p.\,58]) that
$$\sum_{x=0}^{p-1}\l(\f{x^2+ax+b}p\r)=\cases -1&\t{if}\ p\nmid a^2-4b,
\\p-1&\t{if}\ p\mid a^2-4b.\endcases$$
So, for each $j=1,\ldots,p-1$, we have
$$\sum_{i=1}^{p-1}\l(\f{i^2+cij+dj^2}p\r)=\cases-1-(\f dp)&\t{if}\ p\nmid c^2-4d,\\p-1-(\f dp)&\t{if}\ p\mid c^2-4d.
\endcases$$

Now assume that $(\f dp)=1$. Set
$$\lambda=\cases1/2&\t{if}\ p\nmid c^2-4d,
\\1/(2-p)&\t{if}\ p\mid c^2-4d.\endcases$$
Then
$$\lambda\sum_{i=1}^{p-1}\l(\f{i^2+cij+dj^2}p\r)=-1.$$
By adding the last $p-1$ rows multiplied by $\lambda$ to the first row of $[c,d]_p$, the initial term of the resulting determinant becomes
$(p-1)\lambda$ while all other terms in the first row vanish. It follows that
$$[c,d]_p=(p-1)\lambda(c,d)_p$$
which is equivalent to (1.23). This ends the proof. \qed

\heading{3. Proofs of Theorems 1.4-1.5}\endheading

\proclaim{Lemma 3.1} Let $p$ be an odd prime. For any $d\in\Z$ with $(\f dp)=-1$, we have the new congruence
$$\prod_{x=1}^{(p-1)/2}(x^2-d)\eq (-1)^{(p+1)/2}\,2\pmod p.\tag3.1$$
\endproclaim
\Proof. As $1^2,\ldots,(\f{p-1}2)^2$ give all the $\f{p-1}2$ quadratic residues modulo $p$, we have
$$\prod_{x=1}^{(p-1)/2}(z-x^2)\eq z^{(p-1)/2}-1\pmod p.\tag3.2$$
Thus
$$\prod_{x=1}^{(p-1)/2}(d-x^2)\eq d^{(p-1)/2}-1\eq-2\pmod p$$
and hence (3.1) follows. \qed

\medskip
\noindent{\it Proof of Theorem 1.4}. (i)  Set $m=(p-1)/2$. Clearly
$$\bg|\f{(\f{i+j}p)}{i+j}\bg|_{1\ls i,j\ls (p-1)/2}\eq|(i+j)^{m-1}|_{1\ls i,j\ls m}\pmod p.$$
By Lemma 2.1,
$$\align |(i+j)^{m-1}|_{1\ls i,j\ls m}=&\prod_{k=0}^{m-1}\bi{m-1}k\times\prod_{1\ls i<j\ls m}(i-j)(j-i)
\\=&\f{(m-1)!^m}{\prod_{k=0}^{m-1}k!(m-1-k)!}(-1)^{m(m-1)/2}\prod_{1\ls i<j\ls m}(j-i)^2
\\=&(-1)^{m(m-1)/2}(m-1)!^m
\\=&(-1)^{(p-1)(p-3)/8}\f{2^{(p-1)/2}}{(p-1)^{(p-1)/2}}\l(\f{p-1}2!\r)^{(p-1)/2}.
\endalign$$
Since $(\f2p)=(-1)^{(p^2-1)/8}$, by the above we have
$$\bg|\f{(\f{i+j}p)}{i+j}\bg|_{1\ls i,j\ls (p-1)/2}\eq\l(\f{p-1}2!\r)^{(p-1)/2}\pmod p.\tag3.3$$
Recall that
$$\l(\f{p-1}2!\r)^2\eq(-1)^{(p+1)/2}\pmod p$$
by Remark 1.1(a).
In the case $p\eq1\ (\mo\ 4)$, from (3.3) we get
$$\bg|\f{(\f{i+j}p)}{i+j}\bg|_{1\ls i,j\ls (p-1)/2}\eq(-1)^{(p-1)/4}=\l(\f 2p\r)\pmod p.$$
If $p\eq3\ (\mo\ 4)$, then by (3.3) we have
$$\bg|\f{(\f{i+j}p)}{i+j}\bg|_{1\ls i,j\ls (p-1)/2}\eq\f{p-1}2!\l(\f{p-1}2!\r)^{2(p-3)/4}\eq\f{p-1}2!\pmod p.$$
So (1.24) always holds.

(ii) It is known (cf. [K05, (5.5)]) that
$$\bg|\f1{x_i+y_j}\bg|_{1\ls i,j\ls m}=\f{\prod_{1\ls i<j\ls m}(x_j-x_i)(y_j-y_i)}{\prod_{i=1}^m\prod_{j=1}^m(x_i+y_j)}.$$
Taking $m=(p-1)/2$ and $x_i=y_i=i^2$ for $i=1,\ldots,m$, we get
$$\bg|\f1{i^2+j^2}\bg|_{1\ls i,j\ls (p-1)/2}=\f{\prod_{1\ls i<j\ls (p-1)/2}(j^2-i^2)^2}{\prod_{i=1}^{(p-1)/2}\prod_{j=1}^{(p-1)/2}(i^2+j^2)}.\tag3.4$$
Observe that
$$\align\prod_{1\ls i<j\ls(p-1)/2}(j^2-i^2)^2=&\prod_{j=1}^{(p-1)/2}((j-1)!(j+1)\cdots(2j-1))^2
\\=&\prod_{j=1}^{(p-1)/2}\f{(2j-1)!^2}{j^2}=\f{\prod_{j=1}^{(p-1)/2}(2j-1)!(p-2j)!}{((p-1)/2)!^2}
\\=&\f{\bi{p-1}{(p-1)/2}}{(p-1)!}\prod_{j=1}^{(p-1)/2}\f{(p-1)!}{\bi{p-1}{2j-1}}
\\\eq&\f{(-1)^{(p-1)/2}}{-1}\prod_{j=1}^{(p-1)/2}\f{-1}{(-1)^{2j-1}}=1\pmod p
\endalign$$
in view of Wilson's theorem and the congruence $\bi{p-1}k\eq(-1)^k\pmod p$ for $k=0,\ldots,p-1$. Also,
$$\align\prod_{i=1}^{(p-1)/2}\prod_{j=1}^{(p-1)/2}(i^2+j^2)=&\prod_{i=1}^{(p-1)/2}\(\prod_{j=1}^{(p-1)/2}i^2\l(1+\f{j^2}{i^2}\r)\)
\\\eq&\prod_{i=1}^{(p-1)/2}\(i^{p-1}\prod_{x=1}^{(p-1)/2}(1+x^2)\)\pmod p.
\endalign$$
As $-1$ is a quadratic non-residue modulo $p$, applying (3.1) with $d=-1$ we get
$$\prod_{x=1}^{(p-1)/2}(x^2+1)\eq(-1)^{(p+1)/2}2=2\pmod p.$$
Therefore
$$\prod_{i=1}^{(p-1)/2}\prod_{j=1}^{(p-1)/2}(i^2+j^2)\eq\prod_{i=1}^{(p-1)/2}2=2^{(p-1)/2}\eq\l(\f2p\r)\pmod p.$$
So the desired congruence (1.25) follows from (3.4). We are done. \qed

\proclaim{Lemma 3.2} Let $p$ be an odd prime. If $p\eq1\pmod4$, then
$$\l(\f{((p-3)/2)!!}p\r)=(-1)^{|\{0<k<\f p4:\ (\f kp)=-1\}|}.\tag3.5$$
If $p\eq3\pmod4$, then
$$\l(\f{((p-2+(\f 2p))/2)!!}p\r)=(-1)^{\lfloor(p+1)/8\rfloor}.\tag3.6$$
\endproclaim
\Proof. We first handle the case $p\eq1\pmod4$. By Lemma 2.3,
$$\l(\f 2p\r)=\l(\f{((p-1)/2)!}p\r)=\l(\f{((p-1)/2)!!((p-3)/2)!!}p\r).$$
Thus
$$\align\l(\f{((p-3)/2)!!}p\r)=&\l(\f 2p\r)\l(\f{\prod_{0<k<p/4}(2k)}p\r)
\\=&\l(\f 2p\r)^{(p+3)/4}(-1)^{|\{0<k<\f p4:\ (\f kp)=-1\}|}.
\endalign$$
This implies (3.5) since
$$\l(\f 2p\r)^{(p+3)/4}=(-1)^{\f{p+1}2\cdot\f{p-1}4\cdot\f{p+3}4}=1.$$

If $p\eq3\pmod8$, then $\sum_{0<k<p/4}(\f kp)=0$ by [BC],
and hence
$$\l(\f{((p-3)/2)!!}p\r)=\l(\f{\prod_{0<k<p/4}(2k)}p\r)=(-1)^{|\{1\ls k\ls\f {p-3}4:\ (\f kp)=1\}|}=(-1)^{(p-3)/8}$$
since $(\f 2p)=-1$. Similarly, if $p\eq7\pmod8$, then $\sum_{p/4<k<p/2}(\f kp)=0$ by [BC], and hence
$$\align\l(\f{((p-1)/2)!!}p\r)=&(-1)^{|\{1\ls k\ls\f {p-1}2:\ 2\nmid k\ \t{and}\ (\f kp)=-1\}|}
\\=&(-1)^{|\{1\ls k\ls\f {p-1}2:\ (\f kp)=-1\}|-|\{1\ls j\ls \f{p-3}4:\ (\f{2j}p)=-1\}|}
\\=&(-1)^{|\{\f p4<k<\f p2:\ (\f kp)=-1\}|}=(-1)^{\f12|\{k\in\Z:\ \f{p-3}4<k\ls\f{p-1}2\}|}
\\=&(-1)^{(p+1)/8}.
\endalign$$
So (3.6) holds when $p\eq3\pmod4$. \qed

\medskip
\noindent{\it Proof of Theorem 1.5}. Set $m:=(p-1)/2$. Applying Lemma 2.1 with $n=m+1$ and $P(z)=z^m$, we obtain
$$\align W_p\eq&\l|(i^2-m!j)^m\r|_{0\ls i,j\ls m}=\prod_{k=0}^{m}\bi mk\times\prod_{0\ls i<j\ls m}(i^2-j^2)(-m!j-(-m!i))
\\=&\prod_{k=0}^{m}\bi mk\times(m!)^{m(m+1)/2}\prod_{0\ls i<j\ls m}(j-i)^2(i+j)\pmod p.
\endalign$$
Observe that
$$\align\prod_{0\ls i<j\ls m}(i+j)=&\prod_{k=1}^m k^{|\{i\in\Z:\ 0\ls i\ls\lfloor\f{k-1}2\rfloor\}|}
(p-k)^{|\{i\in\Z:\ m+1-k\ls i\ls\lfloor m-\f k2\rfloor\}|}
\\=&\prod_{k=1}^mk^{\lfloor (k+1)/2\rfloor}(p-k)^{\lfloor k/2\rfloor}
\\\eq&(-1)^{\sum_{k=1}^m\lfloor k/2\rfloor}\prod_{k=1}^mk^k
=(-1)^{m-m/(2-\{m\}_2)}\prod_{k=1}^mk^k\pmod p.
\endalign$$
So we have
$$W_p\eq(-1)^{m-m/(2-\{m\}_2)}(m!)^{m(m+1)/2}\prod_{k=0}^m\bi mk
\times\prod_{j=1}^m (j!)^2\times\prod_{k=1}^m k^k\pmod p.\tag3.7$$

{\it Case} 1. $p\eq1\pmod4$ (i.e., $2\mid m$).

In this case,
$$\prod_{k=0}^m\bi mk=\prod_{k=0}^{m/2-1}\bi mk^2\times\bi{m}{m/2}.$$
Write $p=x^2+4y^2$ with $x,y\in\Z$ and $x\eq1\pmod4$. By a result of Gauss (cf. [BEW, (9.0.1)]),
$$\bi{(p-1)/2}{(p-1)/4}\eq2x\pmod p\tag3.8$$
and hence
$$\l(\f{\bi{m}{m/2}}p\r)=\l(\f{2x}p\r)=\l(\f2p\r)\l(\f{|x|}p\r)=\l(\f 2p\r)\l(\f p{|x|}\r)=\l(\f 2p\r)\l(\f{x^2+4y^2}{|x|}\r)=\l(\f 2p\r).$$
Combining this with (3.7), (2.4) and (3.5), we obtain that
$$\align\l(\f{W_p}p\r)=&\l(\f{-1}p\r)^{m/2}\l(\f 2p\r)\l(\f{m!}p\r)^{m(m+1)/2}\l(\f{\prod_{k=1}^mk^k}p\r)
\\=&\l(\f 2p\r)^{1+(p^2-1)/8}\l(\f{(m-1)!!}p\r)=\l(\f{((p-3)/2)!!}p\r)
\\=&(-1)^{|\{0<k<\f p4:\ (\f kp)=-1\}|}.
\endalign$$

{\it Case} 2. $p\eq3\pmod4$ (i.e., $2\nmid m$).

In this case,
$$\prod_{k=0}^m\bi mk=\prod_{k=0}^{(m-1)/2}\bi mk^2.$$
Combining this with (3.7) and (3.6), we see that
$$\align\l(\f{W_p}p\r)=&\l(\f{m!}p\r)^{m(m+1)/2}\l(\f{\prod_{k=1}^mk^k}p\r)
=\l(\f{m!}p\r)^{(p^2-1)/8}\l(\f{m!!}p\r)
\\=&\cases(\f{m!}p)(\f{m!!}p)=(\f{(m-1)!!}p)=(-1)^{\lfloor(p+1)/8\rfloor}&\t{if}\ p\eq3\pmod 8,
\\(\f{m!!}p)=(-1)^{(p+1)/8}&\t{if}\ p\eq7\pmod 8.\endcases
\endalign$$

In view of the above, we have completed the proof of (1.26). \qed

\heading{4. Some open conjectures on determinants}\endheading

By Remark 1.1(a),
$$\f{p-1}2!\eq\pm1\pmod p\quad\t{for any prime}\ p\eq3\pmod 4.$$
In 1961 L. J. Mordell [M61] proved that
for any prime $p>3$ with $p\eq3\ (\mo\ 4)$ we have
$$\f{p-1}2!\eq(-1)^{(h(-p)+1)/2}\pmod p,$$
where $h(-p)$ is the class number of the imaginary quadratic field $\Q(\sqrt{-p})$.

Motivated by Theorem 1.5, we pose the following conjecture which looks quite challenging.

\proclaim{Conjecture 4.1} Let $p$ be an odd prime. Then
$$\bg|\l(\f{i^2-((p-1)/2)!j}p\r)\bg|\Sb 1\ls i,j\ls(p-1)/2\endSb=0\ \ \text{if and only if}\ \ p\eq 3\ (\mo\ 4).\tag4.1$$
When $p\eq3\pmod4$, both
$$(-1)^{\lfloor(p+1)/8\rfloor}\bg|x+\l(\f{i^2-((p-1)/2)!j}p\r)\bg|_{0\ls i,j\ls(p-1)/2}$$
and
$$\f{(-1)^{\lfloor(p+1)/8\rfloor}}x\bg|x+\l(\f{i^2-((p-1)/2)!j}p\r)\bg|_{1\ls i,j\ls(p-1)/2}$$
are positive squares not depending on $x$.
\endproclaim

\Remark\ 4.1. We have verified that (4.1) holds for all odd primes $p<2500$,
and that $(-1)^{\lfloor(p+1)/8\rfloor}W_p$ is a positive square for every prime $p<3000$
with $p\eq3\pmod4$. See [Su13, A226163] for the sequence
$$\bg|\l(\f{i^2-((p_n-1)/2)!j}{p_n}\r)\bg|_{1\ls i,j\ls(p_n-1)/2}\ \ (n=2,3,\ldots),$$
where $p_n$ denotes the $n$th prime.
\medskip

Inspired by Theorem 1.2(iii), we formulate the following conjecture.

\proclaim{Conjecture 4.2} Let $p$ be an odd prime.

{\rm (i)} Let $S^*(1,p)$
denote the determinant obtained from $S(1,p)=|(\f{i^2+j^2}p)|_{1\ls i,j\ls(p-1)/2}$
via replacing the entries $(\f{1^2+j^2}p)\ (j=1,\ldots,(p-1)/2)$ in the first row
by $(\f jp)\ (j=1,\ldots,(p-1)/2)$ respectively. Then $-S_p$ is always a positive square, where
$$S_p:=\cases S^*(1,p)&\t{if}\ p\eq1\pmod4,\\S(1,p)/2^{(p-3)/2}&\t{if}\ p\eq3\pmod4.\endcases\tag4.2$$

{\rm (ii)} If $p=x^2+4y^2$ with $x$ and $y$ positive integers, then for any $d\in\Z$ with $(\f dp)=-1$, the number
$|T(d,p)|/(2^{(p-1)/2}y)$ is a positive square not depending on $d$.
\endproclaim
\Remark\ 4.2. We have verified this
for all odd primes $p<2000$; for example,
$$S_3=S_5=S_7=S_{11}=-1^2,\ S_{13}=-3^2,\ S_{17}=-21^2,\ S_{19}=-2^2,\ S_{23}=-1^2.$$
 After learning part (i) of this conjecture in the case $p\eq3\pmod4$ from the author, H. Cohen and M. Vsemirnov informed the author their supplemental conjecture which states that for any prime $p=x^2+4y^2$ with $x,y\in\Z$ and $x\eq1\pmod4$, the number $S(1,p)/x$ is an integer square.
\medskip

The following conjecture can be viewed as a supplement to Theorem 1.2(iii).

\proclaim{Conjecture 4.3} Let $p$ be an odd prime, and let $c,d\in\Z$ with $p\nmid cd$. Define
$$S_c(d,p):=\bg|\l(\f{i^2+dj^2+c}p\r)\bg|\Sb 1\ls i,j\ls(p-1)/2\endSb.\tag4.3$$
Then
$$\l(\f {S_c(d,p)}p\r)=\cases1&\t{if}\ (\f cp)=1\ \t{and}\ (\f dp)=-1,
\\(\f{-1}p)&\t{if}\ (\f cp)=(\f dp)=-1,
\\(\f{-2}p)&\t{if}\ (\f{-c}p)=(\f dp)=1,
\\(\f{-6}p)&\t{if}\ (\f{-c}p)=-1\ \t{and}\ (\f dp)=1.\endcases\tag4.4$$
\endproclaim
\Remark\ 4.3. We have verified Conjecture 4.3 for all odd primes $p<225$. See [Su13, A228005] for the sequence $S_1(1,p_n)\ (n=2,3,\ldots)$.
Let $p$ be an odd prime and let $b,c,d\in\Z$ with $p\nmid bcd$. It is easy to see that
$$\align\bg|\l(\f{i^2+dj^2+b^2c}p\r)\bg|\Sb 1\ls i,j\ls(p-1)/2\endSb
=&\bg|\l(\f{(bi)^2+d(bj)^2+b^2c}p\r)\bg|\Sb 1\ls i,j\ls(p-1)/2\endSb
\\=&\bg|\l(\f{i^2+dj^2+c}p\r)\bg|\Sb 1\ls i,j\ls(p-1)/2\endSb.
\endalign$$
\medskip

\proclaim{Conjecture 4.4} For any prime $p\eq\pm3\pmod 8$, the number
$$-\l|\l(\f{i(i+1)+j(j+1)}p\r)\r|_{1\ls i,j\ls(p-1)/2}$$
is a quadratic nonresidue modulo $p$.
\endproclaim
\Remark\ 4.4. We have verified this conjecture for $p<1000$.
Let $p$ be an odd prime. By using Lemma 2.1, and Gauss' congruence (3.8) for $p\eq1\pmod4$,
we are able to show that $2|(\f{i(i+1)+j(j+1)}p)|_{0\ls i,j\ls(p-1)/2}$ is always a quadratic residue modulo $p$.
\medskip

\proclaim{Conjecture 4.5} {\rm (i)} For any positive odd integer $n>3$, we have
$$\bg|(i^2+j^2)\l(\f{i^2+j^2}n\r)\bg|_{0\ls i,j\ls (n-1)/2}\eq0\pmod n.\tag4.5$$

{\rm (ii)} Let $p>5$ be a prime with $p\eq1\ (\mo\ 4)$. Then
$$\l(\f{D_p^+}p\r)=\l(\f{D_p^-}p\r)=1,\tag4.6$$
where
$$D_p^+:=\bg|(i+j)\l(\f{i+j}p\r)\bg|_{1\ls i,j\ls(p-1)/2}\ \t{and}\ D_p^-:=\bg|(j-i)\l(\f{j-i}p\r)\bg|_{1\ls i,j\ls(p-1)/2}.$$

{\rm (iii)} For any prime $p>3$ and integer $d\not\eq0\pmod p$, we have
$$\left(\frac{|(i^2+dj^2)^{(p+1)/2}|_{1\ls i,j\ls(p-1)/2}}p\right)
=\cases(\f dp)^{(p-1)/4}&\text{if}\ 4\mid p-1,
\\(\f dp)^{(p+1)/4}(-1)^{(h(-p)-1)/2}&\text{if}\ 4\mid p-3.\endcases$$
\endproclaim
\Remark\ 4.5. We have verified part (i) and parts (ii)-(iii) of Conjecture 4.5
for $n<1200$ and $p<1200$ respectively. It is known that a skew-symmetric determinant of even order with integer entries is always a square (cf. [St90] and [K99]).

\proclaim{Conjecture 4.6} Let $p$ be an odd prime.

{\rm (i)} Define $M_p$ as the matrix obtaining from $[(\frac{i-j}p)]_{0\ls i,j\ls (p-1)/2}$ via replacing all the entries in the first row by $1$. Then
$$\det M_p=\cases(-1)^{(p-1)/4}&\text{if}\ p\equiv 1\pmod4,
\\(-1)^{(h(-p)-1)/2}&\text{if}\ p>3\ \text{and}\ p\equiv3\pmod4.\endcases\tag4.7$$

{\rm (ii)} Let $M_p^+$ be the matrix obtaining from
$[(\frac{i+j}p)]_{0\ls i,j\ls(p-1)/2}$ via replacing all the entries in the first row by $1$. Then
$$\det M_p^+=2^{(p-1)/2}\det M_p.\tag4.8$$

{\rm (iii)} Define $N_p$ as the matrix obtaining from $[(\frac{i-j}p)]_{1\ls i,j\ls (p-1)/2}$ via replacing all the entries in the first row by $1$, and also define $N_p^+$ as the matrix obtaining from $[(\frac{i+j}p)]_{1\ls i,j\ls (p-1)/2}$ via replacing all the entries in the first row by $1$. Then
$$\det N_p=(-1)^{\lfloor (p+3)/4\rfloor}\det M_p\ \ \text{and}\
\ \det N_p^+=-2^{(p-3)/2}\det M_p.\tag4.9$$
\endproclaim
\Remark\ 4.6. We have verified Conjecture 4.6 for all odd primes $p<2000$.
\medskip

Recall that for any $c,d\in\Z$ and positive odd integer $n$,
 $(c,d)_n$ and $[c,d]_n$ are defined by (1.21) and (1.22).

\proclaim{Conjecture 4.7} Let $c$ and $d$ be integers.

{\rm (i)} If $d$ is nonzero, then $[c,d]_p=0$
for infinitely many primes $p$.

{\rm (ii)} If $[c,d]_p$ is nonzero for some odd prime $p$, then its $p$-adic valuation $($i.e., $p$-adic order$)$ must be even.
\endproclaim
\Remark\ 4.7. We have verified Conjecture 4.7(ii) for all odd primes $p<200$. Let $c,d\in\Z$. If $d$ is not a square, then by [IR, pp.\,57-58] there are infinitely many odd primes $p$
with $(\f dp)=-1$ and hence $(c,d)_p=0$ by Theorem 1.3(i). If $p$ is an odd prime with $p\nmid d$ and $(c,d)_p\not=0$, then
$(\f dp)=1$, $[c,d]_p\not=0$ and $\ord_p(c,d)_p=\ord_p[c,d]_p$ by Theorem 1.3, where $\ord_p x$ denotes the $p$-adic order of a nonzero $p$-adic number $x$.
\medskip

\proclaim{Conjecture 4.8} {\rm (i)} For any odd integer $n>3$ we have $(2,3)_n\eq0\pmod{n^2}$. Also, for any odd integer $n>5$ we have
$(6,15)_n\eq0\pmod {n^2}$.

{\rm (ii)} For any positive integer $n\eq3\pmod4$, we have
$$(6,1)_n=[6,1]_n=(3,2)_n=[3,2]_n=0\tag4.10$$
and
$$(4,2)_n=(8,8)_n=(3,3)_n=(21,112)_n=0.\tag4.11$$

{\rm (iii)} For any prime $p\eq1\pmod4$ with $p=x^2+4y^2\ (x,y\in\Z\ \t{and}\ 4\mid x-1)$, we have
$$(6,1)_p=(-1)^{(p-1)/4}xu^2\ \ \t{and}\ \ (3,2)_p=(-1)^{(p-1)/4}xv^2$$
for some $u,v\in\Z$. For each prime $p\eq1\pmod 8$ with $p=x^2+2y^2\ (x,y\in\Z\ \t{and}\ 4\mid x-1)$, we have
$(4,2)_p=(8,8)_p=(-1)^{(p-1)/8}xw^2$ for some $w\in\Z$. For every prime $p\eq1\pmod {12}$ with $p=x^2+3y^2\ (x,y\in\Z\ \t{and}\ 3\mid x-1)$, the number
$(3,3)_p/x$ is a square. Also, for any prime $p\eq1,9,25\pmod{28}$ with $p=x^2+7y^2\ (x,y\in\Z\ \t{and}\ (\f x7)=1)$ the number $(21,112)_p/x$ is a square.

{\rm (iv)} $(10,9)_p=0$ for any prime $p\eq5\pmod{12}$, $[5,5]_p=0$ for each prime $p\eq13,17\pmod{20}$,
and $(8,18)_p=[8,18]_p=0$ for every prime $p\eq13,17\pmod{24}$.
\endproclaim
\Remark\ 4.8. We have verified parts (i) and (ii) of Conjecture 4.8 for $n<1000$, and parts (iii) and (iv)
of Conjecture 4.8 for primes $p<1000$. See [Su13, A225611] for the sequence $(6,1)_{p_n}\ (n=2,3,\ldots)$. For any positive integer $n\eq\pm5\pmod{12}$, we have $(\f 3n)=-1$ and hence $(2,3)_n=(3,3)_n=0$ by Theorem 1.3(i).
Similarly, for any positive integer $n\eq\pm3\pmod{8}$, we have $(\f 2n)=-1$ and hence $(3,2)_n=(4,2)_n=(8,8)_n=0$ by Theorem 1.3(i). Also, for any positive integer $n\eq3,7\pmod{10}$,
we have $(\f 5n)=-1$ and hence $(5,5)_n=0$ by Theorem 1.3(i). For any prime $p\eq5\pmod{12}$, clearly $(10,9)_p=0$
if and only if $[10,9]_p=0$ (by Theorem 1.3(ii)).

\Ack. The author would like to thank Prof. R. Chapman and C. Krattenthaler, and his former students X. Meng, H. Pan and L. Zhao, and the two anonymous referees for their helpful comments.

\enddocument